\documentclass[11pt]{article}

\advance \topmargin by -\headheight
\advance \topmargin by -\headsep
\evensidemargin \oddsidemargin
\usepackage{amssymb,latexsym,amsmath,epsfig,amsthm,tikz} 
\usepackage{ifthen}

\usepackage{color}
\definecolor{darkgreen}{rgb}{0,.4,0.2}
\definecolor{darkagenta}{rgb}{.5,0,.5}
\definecolor{darkred}{rgb}{0.85,0,0}
\definecolor{darkblue}{rgb}{0,0,.6}
\definecolor{lightgray}{gray}{.95}
\definecolor{rgrey}{rgb}{.8,0.4,.4}  
\definecolor{grey}{rgb}{.5,.5,.5}  

\newtheorem{theorem}{Theorem}
\newtheorem{lemma}{Lemma}
\newtheorem{proposition}{Proposition}

\theoremstyle{definition}

\newcommand{\A}{ {\cal A}}

\begin{document}

\begin{center}
\uppercase{\bf Sumsets and fixed points of substitutions}
\vskip 20pt
{\bf F.~Michel Dekking}\\

\bigskip

{\rm DIAM,  Delft University of Technology, Faculty EEMCS, Delft, The Netherlands.}\\

\end{center}

\bigskip

\centerline{\bf Abstract}

\noindent In this paper we introduce a technique to determine the sumset $A+A$, where $A$ is the indicator function of the 0's occurring in a fixed point $x$ of a substitution on the alphabet $\{0,1\}$.

\vskip 30pt

\section{Introduction}

Let $A$ be the sequence given by $A(n)= \lfloor n\alpha\rfloor$ for $n\ge 1$, where $\alpha = (1+\sqrt{5})/2$ is the golden mean. In the recent paper \cite{Kawsumarng}, it was proved that the sumset $A+A$ is equal to the natural numbers with exception of the numbers 1 and 3.
Here the sumset $A+A$ is defined by
$$A+A = \{a+b: a\in A, b\in A\}.$$
It was also shown that if $B$ is defined by $B(n)= \lfloor n\alpha^2\rfloor$, then  the set $B+B$ has a complement that is infinite. The determination of $B+B$ takes more than 8 pages, and the  complicated structure of $B+B$ is described as containing ``some fractal patterns". The goal of the present paper is to elucidate what the source of these fractal patterns is, and at the same time to introduce a vast generalization of these types of sumsets.

\medskip

Interestingly, Shallit in \cite{Shallit-sumset} also positioned the sumset problem of the paper \cite{Kawsumarng} in the combinatorics on words context, but in a completely different way. The focus there is on the Fibonacci representation (also known as Zeckendorf representation) of the natural numbers, and the proofs are computer assisted. See also the paper \cite{Shallit-additive}.

\medskip

The crucial observation that we make is that the sequences $A=(\lfloor n\alpha\rfloor) =(1, 3, 4, 6, 8, 9, 11,\dots)$ and $B=(\lfloor n\alpha^2\rfloor)=(2, 5, 7, 10, 13, 15,\dots)$  give the positions of 0's, respectively the positions of 1's in the infinite Fibonacci word $x_{\rm F}=010010100100\dots$, fixed point of the substitution $0\mapsto 01, \, 1\mapsto 0$ (see, e.g., \cite{Lothaire}). Let $\sigma$ be any substitution on the monoid $\{0,1\}^*$, admitting a fixed point $x=\sigma(x)$. Then, let
$A$ give the positions of 0 in $x$, and let $B$ give the positions of 1 in $x$.

Problem: determine $A+A$, $A+B$ and $B+B$.

\medskip

 The proofs in \cite{Kawsumarng} are purely arithmetical. Our approach, inspired by \cite{DSS-AP}, is to analyse  $A+A$ by studying the product set $A\times A$.
 This amounts to passing from fixed points  of substitutions to fixed points of two-dimensional substitutions.  See, e.g., \cite{Frank-primer} for an overview of the theory of two-dimensional substitutions. Our application of two-dimensional substitutions to the sumset problem is self-contained.

 \medskip  In Section \ref{sec:Fib} we give a very simple proof of the  $A+A$  result of \cite{Kawsumarng} using a
 two-dimensional substitution which has $x_{\rm F}\times x_{\rm F}$ as fixed point. 
A proof of the $B+B$ result from \cite{Kawsumarng}, analogously to the way this is done in  Theorem \ref{th:vN}, is very complex. The numbers of columns on the right side of  the 2D fixed point of the product substitution that one has to use is unbounded for the Fibonacci sunstitution, whereas it is 15 for the von Neumann substitution.

\medskip

In Section \ref{sec:TM} we solve the $A+A$, $A+B$ and $B+B$ problem for the Thue-Morse  word $t=0110100110010110\dots$, fixed point of  the substitution $0\mapsto 01, \, 1\mapsto 10$.
We remark that our result is equivalent to Theorem 2 of the paper \cite{Shallit-additive}, which is proved in a completely different, computer assisted way.

\medskip

In Section \ref{sec:vonNeu} we solve the $A+A$ and $B+B$ problem for the von Neumann word.
The interest of this example lies in the fact that this example is one of a large family for which the method of \cite{Shallit-sumset} does not work (\cite{Shallit-Walnut}).

\medskip

In Section \ref{sec:sq} we present the classical `sum of two squares' problem in our 2D substitution context. This example also illustrates the fact that our technique extends from fixed points  of substitutions to 0-1-valued morphic sequences, i.e.,  fixed points of substitutions on arbitrary alphabets which are mapped to $\{0,1\}^*$ by a letter-to-letter map.

\section{Substitutions and their products} \label{sec:sub}

  A substitution $\sigma$ is a  homomorphism of the monoid $\A^*$ of words over an alphabet $\A$, that is, $\A^*$ consists of concatenations $a_1\dots a_m$ with $a_i$ from $\A$ for $i=1\dots m$, and $\sigma(vw)=\sigma(v)\sigma(w)$ for all $v,w\in \A^*$.
  Since we are interested in the characteristic functions of 0-1-words, we simplify the presentation and take $\A=\{0,1\}$.

  Let $\sigma$ on $\A^*=\{0,1\}^*$ be a substitution given by
  $$\sigma(0)=a_1\dots a_m, \quad \sigma(1)=b_1\dots b_{m'},$$
  for two natural numbers $m$ and $m'$. We then define the direct product substitution $\sigma\!\times\!\sigma$ on the alphabet $\A\times \A$ by
  \begin{align*}
    \sigma\!\times\!\sigma((0,0))_{(k,\ell)}  & =(a_k,a_\ell), \quad \sigma\!\times\!\sigma((0,1))_{(k,\ell')}=(a_k,b_{\ell'}),  \\
    \sigma\!\times\!\sigma((1,0))_{(k',\ell)} & =(b_{k'},a_\ell),   \quad \sigma\!\times\!\sigma((1,1))_{(k',\ell')}=(b_{k'},b_{\ell'}),
  \end{align*}
where $1\le k \le m, 1\le \ell \le m, 1\le k' \le m', 1\le \ell' \le m'.$

\medskip

\noindent {\bf Example}\, Let $\sigma$ be the Fibonacci substitution given by $\sigma(0)=01,\; \sigma(1)=0$ then, abbreviating $(i,j)$ to $ij$ for $i,j=0,1$
the direct product substitution $\sigma\!\times\!\sigma$ is given  by\\
\begin{tikzpicture}[scale=.3,rounded corners=0]
 \def\bricka(#1){\path (#1)+(-1/2,-1/2) coordinate (shP0); \draw (shP0) rectangle +(1,1); \node[circle, inner sep=5.4pt] at (#1)  {\tiny $00$}}
\def\brickb(#1){\path (#1)+(-1/2,-1/2) coordinate (shP0); \draw (shP0) rectangle +(1,1); \node[circle, inner sep=5.4pt] at (#1)  {\tiny $01$}}
\def\brickc(#1){\path (#1)+(-1/2,-1/2) coordinate (shP0); \draw (shP0) rectangle +(1,1); \node[circle, inner sep=5.4pt] at (#1)  {\tiny $10$}}
\def\brickd(#1){\path (#1)+(-1/2,-1/2) coordinate (shP0); \draw (shP0) rectangle +(1,1); \node[circle, inner sep=5.4pt] at (#1)  {\tiny $11$}}
\def\theta(#1,#2){\path (#2) coordinate (P0);
\path (P0)++(0,0) coordinate  (P00); \path (P0)++(1,0) coordinate (P10); \path (P0)++(0,1) coordinate (P01); \path (P0)++(1,1) coordinate (P11);
\ifthenelse{#1=00}{\bricka(P00);  \brickb(P01); \brickc(P10); \brickd(P11)}{};
\ifthenelse{#1=01}{\bricka(P00);  \brickc(P10)}{};
\ifthenelse{#1=10}{\bricka(P00);  \brickb(P01)}{};
\ifthenelse{#1=11}{\bricka(P00)}{};  }; 
\hspace*{3cm}
\path (0,0) coordinate (a);   \bricka(a);  \node at (2,0) {$\mapsto$};  \path (4,0) coordinate (t00);  \theta(00,t00);
\path (8,0) coordinate (b);   \brickb(b);  \node at (10,0) {$\mapsto$};  \path (12,0) coordinate (t01);  \theta(01,t01);
\path (16,0) coordinate (c);  \brickc(c);  \node at (18,0) {$\mapsto$};  \path (20,0) coordinate (t10); \theta(10,t10);
\path (23,0) coordinate (d);  \brickd(d);  \node at (25,0) {$\mapsto$}; \path (27,0) coordinate (t11); \theta(11,t11);
\end{tikzpicture}

\bigskip

Suppose that $\sigma(0)$ has prefix 0. Then $\sigma^n(0)$ converges to a fixed point $x$ of $\sigma$ as $n\rightarrow\infty$, and so
$[\sigma\!\times\!\sigma]^n((0,0))$ converges to the 2D fixed point $x\times x$ of $\sigma\!\times\!\sigma$.

 Thus the element of $x\times x$  with index $(k,\ell)$  contains the symbol $(0,0)$ if and only if $x_k=0, x_\ell=0$ if and only if  $k\in A$ and $\ell\in A$.\\
Let the \emph{diagonal words} $d_n$ be defined for $n\ge 2$ by
\begin{equation}\label{eq:diag}
d_n =\{[x\times x]_{(k,\ell)}: k+\ell =n, \, k\ge 1, \, l\ge 1\}.
\end{equation}
Then $n\in A+A$ if and only if $d_n$ contains a symbol $(0,0)$. So we have obtained the following theorem.

\begin{theorem}\label{th:diag}
Let $\sigma$ be a substitution on $\{0,1\}$, such that $\sigma(0)$ has prefix 0, and let $x$ be the fixed point of $\sigma$ with prefix $0$. Let $A$ be the sequence of positions of $0$ in $x$. Let $d_n$ be the diagonal words occurring in $x\times x$, for $n\ge 2$.
Then $A+A= \{n\ge 2:  \: the\: symbol\:  (0,0)\:  occurs\: in\: d_n\}.$
\end{theorem}

Note that we also have that if $B$ is the sequence of positions of $1$ in $x$, then  $B+B= \{n\ge 2:  \: the\: symbol\:  (1,1)\:  occurs\: in\: d_n\}$, and $A+B= \{n\ge 2:  \: the\: symbol\:  (0,1)\:or\: (1,0)\:  occurs\: in\: d_n\}$.

\section{The Fibonacci word} \label{sec:Fib}

The Fibonacci word $x_{\rm F}$ is the unique fixed point of the substitution $\sigma:\, 0\mapsto 01, \, 1\mapsto 0$.
Let $\varphi:=\sigma\!\times\!\sigma$ be the direct product of $\sigma$.

It is convenient to code \: $a:=(0,0), b:=(0,1), c:=(1,0), d:=(1,1)$, and to assign colors to these four letters so that one obtains more easily an idea of the structure of $\varphi^n(a)$. The direct product substitution $\varphi$ on this new alphabet is given by

\begin{tikzpicture}[scale=.3,rounded corners=0]
\def\bricka(#1){\path (#1)+(-1/2,-1/2) coordinate (shP0);
                \draw [fill=green!30!white, draw=gray,thick]  (shP0) rectangle +(1,1); \node[circle, inner sep=5.4pt] at (#1)  {\tiny $a$}}
\def\brickb(#1){\path (#1)+(-1/2,-1/2) coordinate (shP0);
                \draw [fill=yellow!30!white, draw=gray,thick]  (shP0) rectangle +(1,1); \node[circle, inner sep=5.4pt] at (#1)  {\tiny $b$}}
\def\brickc(#1){\path (#1)+(-1/2,-1/2) coordinate (shP0);
                \draw [fill=blue!30!white, draw=gray,thick]  (shP0) rectangle +(1,1); \node[circle, inner sep=5.4pt] at (#1)  {\tiny $c$}}
\def\brickd(#1){\path (#1)+(-1/2,-1/2) coordinate (shP0);
                \draw [fill=red!30!white, draw=gray,thick]  (shP0) rectangle +(1,1); \node[circle, inner sep=5.4pt] at (#1)  {\tiny $d$}}
\def\theta(#1,#2){\path (#2) coordinate (P0);
\path (P0)++(0,0) coordinate  (P00); \path (P0)++(1,0) coordinate (P10); \path (P0)++(0,1) coordinate (P01); \path (P0)++(1,1) coordinate (P11);
\ifthenelse{#1=00}{\bricka(P00);  \brickb(P01); \brickc(P10); \brickd(P11)}{};
\ifthenelse{#1=01}{\bricka(P00);  \brickc(P10)}{};
\ifthenelse{#1=10}{\bricka(P00);  \brickb(P01)}{};
\ifthenelse{#1=11}{\bricka(P00)}{};  }; 
\hspace*{3cm}
\path (0,0) coordinate (a);   \bricka(a);  \node at (2,0) {$\mapsto$};  \path (4,0) coordinate (t00);  \theta(00,t00);
\path (8,0) coordinate (b);   \brickb(b);  \node at (10,0) {$\mapsto$};  \path (12,0) coordinate (t01);  \theta(01,t01);
\path (16,0) coordinate (c);  \brickc(c);  \node at (18,0) {$\mapsto$};  \path (20,0) coordinate (t10); \theta(10,t10);
\path (23,0) coordinate (d);  \brickd(d);  \node at (25,0) {$\mapsto$}; \path (27,0) coordinate (t11); \theta(11,t11);
\end{tikzpicture}

Here are some examples of 2D words obtained by  iterating  $\varphi$.

\smallskip

\begin{tikzpicture}[scale=.3,rounded corners=0]
\def\bricka(#1){\path (#1)+(-1/2,-1/2) coordinate (shP0);
                \draw [fill=green!30!white, draw=gray,thick]  (shP0) rectangle +(1,1); \node[circle, inner sep=5.4pt] at (#1)  {\tiny $a$}}
\def\brickb(#1){\path (#1)+(-1/2,-1/2) coordinate (shP0);
                \draw [fill=yellow!30!white, draw=gray,thick]  (shP0) rectangle +(1,1); \node[circle, inner sep=5.4pt] at (#1)  {\tiny $b$}}
\def\brickc(#1){\path (#1)+(-1/2,-1/2) coordinate (shP0);
                \draw [fill=blue!30!white, draw=gray,thick]  (shP0) rectangle +(1,1); \node[circle, inner sep=5.4pt] at (#1)  {\tiny $c$}}
\def\brickd(#1){\path (#1)+(-1/2,-1/2) coordinate (shP0);
                \draw [fill=red!30!white, draw=gray,thick]  (shP0) rectangle +(1,1); \node[circle, inner sep=5.4pt] at (#1)  {\tiny $d$}}
\def\theta(#1,#2){\path (#2) coordinate (P0);
\path (P0)++(0,0) coordinate  (P00); \path (P0)++(1,0) coordinate (P10); \path (P0)++(0,1) coordinate (P01); \path (P0)++(1,1) coordinate (P11);
\ifthenelse{#1=00}{\bricka(P00);  \brickb(P01); \brickc(P10); \brickd(P11)}{};
\ifthenelse{#1=01}{\bricka(P00);  \brickc(P10)}{};
\ifthenelse{#1=10}{\bricka(P00);  \brickb(P01)}{};
\ifthenelse{#1=11}{\bricka(P00)}{};  }; 
\def\thetatwX(#1,#2){\path (#2) coordinate (Q0);
\path (Q0)++(0,0) coordinate  (Q00); \path (Q0)++(2,0) coordinate (Q10); \path (Q0)++(0,2) coordinate (Q01); \path (Q0)++(2,2) coordinate (Q11);
\ifthenelse{#1=00}{\theta(00,Q00);  \theta(01,Q01); \theta(10,Q10); \theta(11,Q11)}{};
\ifthenelse{#1=01}{\theta(00,Q00);  \theta(10,Q10)}{};
\ifthenelse{#1=10}{\theta(00,Q00);  \theta(01,Q01)}{};
\ifthenelse{#1=11}{\theta(00,Q00)}{};  }; 
\def\thetafourX(#1,#2){\path (#2) coordinate (R0);
\path (R0)++(0,0) coordinate  (R00); \path (R0)++(3,0) coordinate (R10); \path (R0)++(0,3) coordinate (R01); \path (R0)++(3,3) coordinate (R11);
\ifthenelse{#1=00}{\thetatwX(00,R00);  \thetatwX(01,R01); \thetatwX(10,R10); \thetatwX(11,R11)}{};
\ifthenelse{#1=01}{\thetatwX(00,R00);  \thetatwX(10,R10)}{};
\ifthenelse{#1=10}{\thetatwX(00,R00);  \thetatwX(01,R01)}{};
\ifthenelse{#1=11}{\thetatwX(00,R00)}{};  }; 
\def\thetaeightX(#1,#2){\path (#2) coordinate (S0);
\path (S0)++(0,0) coordinate  (S00); \path (S0)++(5,0) coordinate (S10); \path (S0)++(0,5) coordinate (S01); \path (S0)++(5,5) coordinate (S11);
\ifthenelse{#1=00}{\thetafourX(00,S00);  \thetafourX(01,S01); \thetafourX(10,S10); \thetafourX(11,S11)}{};
\ifthenelse{#1=01}{\thetafourX(00,S00);  \thetafourX(10,S10)}{};
\ifthenelse{#1=10}{\thetafourX(00,S00);  \thetafourX(01,S10)}{};
\ifthenelse{#1=11}{\thetafourX(00,S00)}{};  }; 
\thetatwX(00,(0,0)); \thetafourX(00,(9,0)); \thetaeightX(00,(20,0)); \thetatwX(01,(34,0)); \thetafourX(01,(43,0));
\node[circle, inner sep=1.4pt] at (-2.5,0) {\footnotesize $\varphi^2(a)\!:$} {};
\node[circle, inner sep=1.4pt] at (6.5,0) {\footnotesize $\varphi^3(a)\!:$} {};
\node[circle, inner sep=1.4pt] at (17.5,0) {\footnotesize $\varphi^4(a)\!:$} {};
\node[circle, inner sep=1.4pt] at (31.5,0) {\footnotesize $\varphi^2(b)\!:$} {};
\node[circle, inner sep=1.4pt] at (40.5,0) {\footnotesize $\varphi^3(b)\!:$} {};
\end{tikzpicture}

Note that
\begin{equation}\label{eq:diag}
d_2=a, \:d_3=bc, \:d_4=ada, \:d_5=acba,\: d_6=bcabc.
\end{equation}

We are now in a position to give a completely different proof of Theorem 3.1. of \cite{Kawsumarng}.

\begin{theorem} {\bf \rm (\cite{Kawsumarng})}
Let $A=1, 3, 4, 6, 8, 9, 11, 12, 14, 16,\dots$ be the  sequence $([n\alpha])$, where $\alpha$ is the golden mean. Then
$$A+A=\mathbb{N}\setminus\{1,3\}.$$
\end{theorem}

\noindent\emph{Proof:} Recall from the introduction that the sequence $([n\alpha])$ is equal to the sequence of positions of 0 in the infinite Fibonacci word $x_{\bf F}$.

We use Theorem \ref{th:diag}. Trivially, 1 is not in $A+A$, and we see that $a$ does not occur in $d_3$. Further we see that $a$ does occur in $d_2$, $d_4$ and $d_5$. So it remains to show that all $d_n$ contain a symbol  $a$ for $n\ge 6$.
To see this, note that the first fife columns on the left border of $x_{\rm F}\times x_{\rm F}$ are a concatenation of the 2D words $\varphi^3(a)$ and $\varphi^3(b)$ for all $n$.
So these five columns are composed of blocks $\varphi^3(a)$ on top of $\varphi^3(a)$, $\varphi^3(a)$ on top of $\varphi^3(b)$, and $\varphi^3(b)$ on top of $\varphi^3(a)$:
\begin{center}
\begin{tikzpicture}[scale=.3,rounded corners=0]
\def\bricka(#1){\path (#1)+(-1/2,-1/2) coordinate (shP0);
                \draw [fill=green!30!white, draw=black]  (shP0) rectangle +(1,1); \node[circle, inner sep=5.4pt] at (#1)  {\tiny $a$}}
\def\brickb(#1){\path (#1)+(-1/2,-1/2) coordinate (shP0);
                \draw [fill=yellow!30!white, draw=black]  (shP0) rectangle +(1,1); \node[circle, inner sep=5.4pt] at (#1)  {\tiny $b$}}
\def\brickc(#1){\path (#1)+(-1/2,-1/2) coordinate (shP0);
                \draw [fill=blue!30!white, draw=black]  (shP0) rectangle +(1,1); \node[circle, inner sep=5.4pt] at (#1)  {\tiny $c$}}
\def\brickd(#1){\path (#1)+(-1/2,-1/2) coordinate (shP0);
                \draw [fill=red!30!white, draw=black]  (shP0) rectangle +(1,1); \node[circle, inner sep=5.4pt] at (#1)  {\tiny $d$}}
\def\theta(#1,#2){\path (#2) coordinate (P0);
\path (P0)++(0,0) coordinate  (P00); \path (P0)++(1,0) coordinate (P10); \path (P0)++(0,1) coordinate (P01); \path (P0)++(1,1) coordinate (P11);
\ifthenelse{#1=00}{\bricka(P00);  \brickb(P01); \brickc(P10); \brickd(P11)}{};
\ifthenelse{#1=01}{\bricka(P00);  \brickc(P10)}{};
\ifthenelse{#1=10}{\bricka(P00);  \brickb(P01)}{};
\ifthenelse{#1=11}{\bricka(P00)}{}; }; 
\def\thetatwX(#1,#2){\path (#2) coordinate (Q0);
\path (Q0)++(0,0) coordinate  (Q00); \path (Q0)++(2,0) coordinate (Q10); \path (Q0)++(0,2) coordinate (Q01); \path (Q0)++(2,2) coordinate (Q11);
\ifthenelse{#1=00}{\theta(00,Q00);  \theta(01,Q01); \theta(10,Q10); \theta(11,Q11)}{};
\ifthenelse{#1=01}{\theta(00,Q00);  \theta(10,Q10)}{};
\ifthenelse{#1=10}{\theta(00,Q00);  \theta(01,Q01)}{};
\ifthenelse{#1=11}{\theta(00,Q00)}{};  }; 
\def\thetafourX(#1,#2){\path (#2) coordinate (R0);
\path (R0)++(0,0) coordinate  (R00); \path (R0)++(3,0) coordinate (R10); \path (R0)++(0,3) coordinate (R01); \path (R0)++(3,3) coordinate (R11);
\ifthenelse{#1=00}{\thetatwX(00,R00);  \thetatwX(01,R01); \thetatwX(10,R10); \thetatwX(11,R11)}{};
\ifthenelse{#1=01}{\thetatwX(00,R00);  \thetatwX(10,R10)}{};
\ifthenelse{#1=10}{\thetatwX(00,R00);  \thetatwX(01,R01)}{};
\ifthenelse{#1=11}{\thetatwX(00,R00)}{}; }; 
\thetafourX(00,(-20,60));\thetafourX(00,(-20,55));
\thetafourX(00,(-10,60));\thetafourX(01,(-10,57));
\thetafourX(01,(0,62));\thetafourX(00,(0,57));
\node[circle, inner sep=2.4pt] at (-23,62)  {\footnotesize $\varphi^3(a):$}; \node[circle, inner sep=2.4pt] at (-23,58)  {\footnotesize $\varphi^3(a):$};
\node[circle, inner sep=2.4pt] at (-13,62)  {\footnotesize $\varphi^3(a):$}; \node[circle, inner sep=2.4pt] at (-13,58)  {\footnotesize $\varphi^3(b):$};
\node[circle, inner sep=2.4pt] at (-3,63)  {\footnotesize $\varphi^3(b):$}; \node[circle, inner sep=2.4pt] at (-3,59)  {\footnotesize $\varphi^3(a):$};
\path (-20,60)++(-1/2,-1/2) coordinate (left); \path (left)++(5,0) coordinate (right); \draw[draw=gray,very thick] (left)--(right);
\path (-10,60)++(-1/2,-1/2) coordinate (left); \path (left)++(5,0) coordinate (right); \draw[draw=gray,very thick] (left)--(right);
\path (-0,62)++(-1/2,-1/2) coordinate (left); \path (left)++(5,0) coordinate (right); \draw[draw=gray,very thick] (left)--(right);
\end{tikzpicture}
\end{center}

In all three cases, the diagonal words starting at the left border in the top $\varphi^3(a)$ or $\varphi^3(b)$ will cross a square with a symbol $a$ in one of the four left most columns, which finishes the proof. \hfill $\Box$.

\medskip

\noindent {\bf Remark}\, The proof shows that if $n=a+a'$, with $a,a'$ from $A$, then  $a$ can always be chosen from the set $\{1,3,4\}$.

\section{The Thue-Morse word} \label{sec:TM}

The Thue-Morse word $t$ is the fixed point of the substitution $\theta:\,0\mapsto 01, 1\mapsto 10$. Although it is general practice (as in \cite{Kawsumarng}) to index Beatty sequences starting from $n=1$, and similarly for fixed points of substitutions, this is \emph{not} the case for the Thue-Morse sequence.
The Thue-Morse word $t=t_0t_1\dots=01101001\dots$ is indexed starting from $n=0$. Thus  $A=0,3,5,6,\dots$ gives the positions of 0 in $t$, and $B=1,2,4,7,\dots$ the positions of 1 in $t$.

Let $\overline{0}=1, \overline{1}=0$ be the symmetry operator on $\{0,1\}^*$. Note that $\theta$ is symmetric, i.e., $\overline{\theta(i)}=\theta(\overline{i})$ for $i=0,1$.

The direct product substitution of $\theta$ is the 2D substitution $\tau$  defined on the symbols $(i,j)$ for $i,j=0,1$ (written as $ij$) by\\[-.9cm]
\begin{align*}\label{eq:tau}
                   & i\overline{j}  \;\,  \overline{i}\overline{j}   \; &     \\[-.1cm]
\tau:\;ij\mapsto\; &  ij            \;\,  \overline{i}j.
\end{align*}

 \smallskip

When we code  \: $  a:=00, b:=01, c:=10, d:=11$, and  color the squares of the symbols $a,b,c$ and $d$, then the substitution $\tau$ takes the form

\smallskip

\begin{tikzpicture}[scale=.4,rounded corners=0]
\def\brickZZ(#1){\path (#1)+(-1/2,-1/2) coordinate (shP0);
                \draw [fill=green!30!white, draw=gray,thick]  (shP0) rectangle +(1,1); \node[circle, inner sep=5.4pt] at (#1)  {\small $a$}}
\def\brickZO(#1){\path (#1)+(-1/2,-1/2) coordinate (shP0);
                \draw [fill=yellow!30!white, draw=gray,thick]  (shP0) rectangle +(1,1); \node[circle, inner sep=5.4pt] at (#1)  {\small $b$}}
\def\brickOZ(#1){\path (#1)+(-1/2,-1/2) coordinate (shP0);
                \draw [fill=blue!30!white, draw=gray,thick]  (shP0) rectangle +(1,1); \node[circle, inner sep=5.4pt] at (#1)  {\small $c$}}
\def\brickOO(#1){\path (#1)+(-1/2,-1/2) coordinate (shP0);
                \draw [fill=red!30!white, draw=gray,thick]  (shP0) rectangle +(1,1); \node[circle, inner sep=5.4pt] at (#1)  {\small $d$}}
\def\theta(#1,#2){\path (#2) coordinate (P0);
\path (P0)++(0,0) coordinate  (P00); \path (P0)++(1,0) coordinate (P10); \path (P0)++(0,1) coordinate (P01); \path (P0)++(1,1) coordinate (P11);
\ifthenelse{#1=00}{\brickZZ(P00);  \brickZO(P01); \brickOZ(P10); \brickOO(P11)}{};
\ifthenelse{#1=01}{\brickZO(P00);  \brickZZ(P01); \brickOO(P10); \brickOZ(P11)}{};
\ifthenelse{#1=10}{\brickOZ(P00);  \brickOO(P01); \brickZZ(P10); \brickZO(P11)}{};
\ifthenelse{#1=11}{\brickOO(P00);  \brickOZ(P01); \brickZO(P10); \brickZZ(P11)}{};
}; 
\path (0,0) coordinate (ZZ);   \brickZZ(ZZ);  \node at (1.5,0) {$\mapsto$};  \path (3,0) coordinate (t00);  \theta(00,t00);
\path (7,0) coordinate (ZO);   \brickZO(ZO);  \node at (8.5,0) {$\mapsto$};  \path (10,0) coordinate (t01);  \theta(01,t01);
\path (14,0) coordinate (OZ);  \brickOZ(OZ);  \node at (15.5,0) {$\mapsto$};  \path (17,0) coordinate (t10); \theta(10,t10);
\path (21,0) coordinate (OO);  \brickOO(OO);  \node at (22.5,0) {$\mapsto$}; \path (24,0) coordinate (t11); \theta(11,t11);
\node[circle, inner sep=1.4pt] at (5,0)  {,}; \node[circle, inner sep=1.4pt] at (12,0)  {,};
\node[circle, inner sep=1.4pt] at (19,0)  {,}; \node[circle, inner sep=1.4pt] at (26,0)  {.};
\end{tikzpicture}

Since we start $t$ at index 0, we have to reconsider the definition of the diagonal words. It turns out that it is convenient to index these by their lengths. So $d_1=a, d_2=bc$, etc.
The iterate $\tau^8(a)$ with the diagonal words $d_1,d_2,\dots,d_{16}$ indicated by lines, is given by

\begin{center}
\begin{tikzpicture}[scale=.4,rounded corners=0]
\def\brickZZ(#1){\path (#1)+(-1/2,-1/2) coordinate (shP0);
                \draw [fill=green!30!white, draw=gray,thick]  (shP0) rectangle +(1,1); \node[circle, inner sep=5.4pt] at (#1)  {\small $a$}}
\def\brickZO(#1){\path (#1)+(-1/2,-1/2) coordinate (shP0);
                \draw [fill=yellow!30!white, draw=gray,thick]  (shP0) rectangle +(1,1); \node[circle, inner sep=5.4pt] at (#1)  {\small $b$}}
\def\brickOZ(#1){\path (#1)+(-1/2,-1/2) coordinate (shP0);
                \draw [fill=blue!30!white, draw=gray,thick]  (shP0) rectangle +(1,1); \node[circle, inner sep=5.4pt] at (#1)  {\small $c$}}
\def\brickOO(#1){\path (#1)+(-1/2,-1/2) coordinate (shP0);
                \draw [fill=red!30!white, draw=gray,thick]  (shP0) rectangle +(1,1); \node[circle, inner sep=5.4pt] at (#1)  {\small $d$}}

\def\theta(#1,#2){\path (#2) coordinate (P0);
\path (P0)++(0,0) coordinate  (P00); \path (P0)++(1,0) coordinate (P10); \path (P0)++(0,1) coordinate (P01); \path (P0)++(1,1) coordinate (P11);
\ifthenelse{#1=00}{\brickZZ(P00);  \brickZO(P01); \brickOZ(P10); \brickOO(P11)}{};
\ifthenelse{#1=01}{\brickZO(P00);  \brickZZ(P01); \brickOO(P10); \brickOZ(P11)}{};
\ifthenelse{#1=10}{\brickOZ(P00);  \brickOO(P01); \brickZZ(P10); \brickZO(P11)}{};
\ifthenelse{#1=11}{\brickOO(P00);  \brickOZ(P01); \brickZO(P10); \brickZZ(P11)}{};
}; 

\def\thetatwX(#1,#2){\path (#2) coordinate (Q0);
\path (Q0)++(0,0) coordinate  (Q00); \path (Q0)++(2,0) coordinate (Q10); \path (Q0)++(0,2) coordinate (Q01); \path (Q0)++(2,2) coordinate (Q11);
\ifthenelse{#1=00}{\theta(00,Q00);  \theta(01,Q01); \theta(10,Q10); \theta(11,Q11)}{};
\ifthenelse{#1=01}{\theta(01,Q00);  \theta(00,Q01); \theta(11,Q10); \theta(10,Q11)}{};
\ifthenelse{#1=10}{\theta(10,Q00);  \theta(11,Q01); \theta(00,Q10); \theta(01,Q11)}{};
\ifthenelse{#1=11}{\theta(11,Q00);  \theta(10,Q01); \theta(01,Q10); \theta(00,Q11)}{};
}; 

\def\thetafourX(#1,#2){\path (#2) coordinate (R0);
\path (R0)++(0,0) coordinate  (R00); \path (R0)++(4,0) coordinate (R10); \path (R0)++(0,4) coordinate (R01); \path (R0)++(4,4) coordinate (R11);
\ifthenelse{#1=00}{\thetatwX(00,R00);  \thetatwX(01,R01); \thetatwX(10,R10); \thetatwX(11,R11)}{};
\ifthenelse{#1=01}{\thetatwX(01,R00);  \thetatwX(00,R01); \thetatwX(11,R10); \thetatwX(10,R11)}{};
\ifthenelse{#1=10}{\thetatwX(10,R00);  \thetatwX(11,R01); \thetatwX(00,R10); \thetatwX(01,R11)}{};
\ifthenelse{#1=11}{\thetatwX(11,R00);  \thetatwX(10,R01); \thetatwX(01,R10); \thetatwX(00,R11)}{};
}; 

\def\thetaeightX(#1,#2){\path (#2) coordinate (S0);
\path (S0)++(0,0) coordinate  (S00); \path (S0)++(8,0) coordinate (S10); \path (S0)++(0,8) coordinate (S01); \path (S0)++(8,8) coordinate (S11);
\ifthenelse{#1=00}{\thetafourX(00,S00);  \thetafourX(01,S01); \thetafourX(10,S10); \thetafourX(11,S11)}{};
\ifthenelse{#1=01}{\thetafourX(01,S00);  \thetafourX(00,S01); \thetafourX(11,S10); \thetafourX(10,S11)}{};
\ifthenelse{#1=10}{\thetafourX(10,S00);  \thetafourX(11,S01); \thetafourX(00,S10); \thetafourX(01,S11)}{};
\ifthenelse{#1=11}{\thetafourX(11,S00);  \thetafourX(10,S01); \thetafourX(01,S10); \thetafourX(00,S11)}{};
}; 
\thetaeightX(00,(-20,30));
\def\row(#1){
\path (#1) coordinate (start);
\foreach \i in {1,2,...,11} {\path (start)++(12,0) coordinate (start);};};
\def\diag(#1){
\path (#1) coordinate (start);
\foreach \i in {1,2,...,16} {\path (start)++(0,\i) coordinate (up); \path (start)++(\i,0) coordinate (down); \draw[gray] (up)--(down);} };
\diag((-20.5,29.5));
\end{tikzpicture}
 \end{center}

%
%
%
%
%

\smallskip

Our goal is to describe the diagonal words $d_n$, for $n=1,2\dots$. We see that in particular
\begin{equation}\label{eq:diagd15}
d_{1}=a, \:d_2=bc, \:d_3=bdc, \:d_4=adda.
\end{equation}

\begin{proposition}\label{prop:diag}
Let $\sigma$ be the morphism given by $\sigma(ij) = i\overline{j}, \overline{i}j$ for $i,j=0,1$. Then
  $$d_{2n}=\sigma(d_n), \quad  for\; n=1,2,\dots.$$
Let $\beta$ be the 2-to-1 morphism given by $\beta(ij,i'j')=ij'$ for $i,j,i',j'=0,1$. Then
$$d_{2n+1}=\beta(d_{2n+2}), \quad  for\; n=1,2,\dots.$$
\end{proposition}

\emph{Proof:} Let the lines $\Delta_n:=\{(x,y):x+y=n-1\}$ correspond to the diagonal words $d_n$. The line $\Delta_{2n}$ cuts through the diagonals of the $2\times 2$-blocks $\tau(a),\dots,\tau(d)$, and these $2\times 2$-blocks are images under $\tau$ of the symbols $a,\dots,d$ on the line $\Delta_n$. Thus  $d_{2n}=\sigma(d_n)$ now follows from the observation that the $i\overline{j}, \overline{i}j$ are exactly the symbols on the diagonals of the $2\times 2$ blocks $\tau(ij)$.

\smallskip

For the second statement, we look at the successive $2\times 2$ blocks $\tau(a),\dots,\tau(d)$  cut by  the line $\Delta_{2n+2}$. The line $\Delta_{2n+1}$ cuts through the lower left squares of these blocks. This implies that the symbols at the odd indices of the word $d_{2n+1}$  can be read off directly from the successive pairs of symbols of  $d_{2n+2}$ by mapping the diagonals of the $2\times 2$ $\tau$-blocks to the symbols at the lower left corner of these blocks.
This gives \, $\beta(i\overline{j},\,\overline{i}j)=   ij$, which is an instance of $\beta(ij,i'j')=ij'$.

\smallskip

For the symbols at the even indices of the word $d_{2n+1}$ we have to do more work: we have to determine the symbol at the upper right corner of the $\tau$-block directly below a pair of $\tau$-blocks from the pair of symbols situated at lower left and the upper right corner of these blocks. The situation is as in the following drawing:

\medskip

\begin{center}
  \begin{tikzpicture}[scale=.5,rounded corners=0]
\def\brickZZ(#1){\path (#1)+(-1/2,-1/2) coordinate (shP0);
                \draw [draw=black]  (shP0) rectangle +(1,1); \node[circle, inner sep=5.4pt] at (#1)  {\small $a$}}
\def\brickZO(#1){\path (#1)+(-1/2,-1/2) coordinate (shP0);
                \draw [draw=black]  (shP0) rectangle +(1,1); \node[circle, inner sep=5.4pt] at (#1)  {\small $b$}}
\def\brickOZ(#1){\path (#1)+(-1/2,-1/2) coordinate (shP0);
                \draw [draw=black]  (shP0) rectangle +(1,1); \node[circle, inner sep=5.4pt] at (#1)  {\small $c$}}
\def\brickOO(#1){\path (#1)+(-1/2,-1/2) coordinate (shP0);
                \draw [draw=black]  (shP0) rectangle +(1,1); \node[circle, inner sep=5.4pt] at (#1)  {\small $d$}}
\def\brickZZempty(#1){\path (#1)+(-1/2,-1/2) coordinate (shP0); \draw [draw=black]  (shP0) rectangle +(1,1); }
\def\brickZOempty(#1){\path (#1)+(-1/2,-1/2) coordinate (shP0); \draw [draw=black]  (shP0) rectangle +(1,1); }
\def\brickOZempty(#1){\path (#1)+(-1/2,-1/2) coordinate (shP0); \draw [draw=black]  (shP0) rectangle +(1,1); }
\def\brickOOempty(#1){\path (#1)+(-1/2,-1/2) coordinate (shP0); \draw [draw=black]  (shP0) rectangle +(1,1); }
\def\theta(#1,#2){\path (#2) coordinate (P0);
\path (P0)++(0,0) coordinate  (P00); \path (P0)++(1,0) coordinate (P10); \path (P0)++(0,1) coordinate (P01); \path (P0)++(1,1) coordinate (P11);
\ifthenelse{#1=00}{\brickZZ(P00);  \brickZO(P01); \brickOZ(P10); \brickOO(P11)}{};
\ifthenelse{#1=01}{\brickZO(P00);  \brickZZ(P01); \brickOO(P10); \brickOZ(P11)}{};
\ifthenelse{#1=10}{\brickOZ(P00);  \brickOO(P01); \brickZZ(P10); \brickZO(P11)}{};
\ifthenelse{#1=11}{\brickOO(P00);  \brickOZ(P01); \brickZO(P10); \brickZZ(P11)}{};
}; 
\def\thetaempty(#1,#2){\path (#2) coordinate (P0);
\path (P0)++(0,0) coordinate  (P00); \path (P0)++(1,0) coordinate (P10); \path (P0)++(0,1) coordinate (P01); \path (P0)++(1,1) coordinate (P11);
\ifthenelse{#1=10}{\brickOZempty(P00);  \brickOOempty(P01); \brickZZempty(P10); \brickZOempty(P11)}{};
}; 


\path (6,4) coordinate (OZ);  \thetaempty(10,OZ);  \path (8,2) coordinate (ZZ); \thetaempty(10,ZZ);  \path (6,2) coordinate (OZ); \thetaempty(10,OZ);
\node[circle, inner sep=1.4pt] at (8,3)  {\footnotesize $i'\overline{j'}$};
\node[circle, inner sep=5.4pt] at (7,4)  {\footnotesize $\overline{i}j$};

\path (12,4) coordinate (OZ);  \thetaempty(10,OZ);  \path (14,2) coordinate (ZZ); \thetaempty(10,ZZ);  \path (12,2) coordinate (OZ); \thetaempty(10,OZ);
\node[circle, inner sep=1.4pt] at (14,3)  {\footnotesize $i'\overline{j'}$ };
\node[circle, inner sep=5.4pt] at (13,4)  {\footnotesize $\overline{i}j$};
\node[circle, inner sep=5.4pt] at (13,3)  {\footnotesize $\overline{i}\,\overline{j'}$};
\node[circle, inner sep=5.4pt] at (10.5,3)  {\footnotesize $\Rightarrow$};
\end{tikzpicture}
\end{center}

\smallskip

This gives: $\beta(\overline{i}j, \footnotesize i'\overline{j'}) = \overline{i}\,\overline{j'}$,
which happens to be a version of $\beta(ij,i'j')=ij'$. \qquad$\Box$

\medskip

\begin{theorem}\label{th:TM}
Let $A=0,3,5,6,9,10,\dots$ give the positions of $0$'s in Thue Morse sequence, and $B=1,2,4,6,7,8,\dots$ the positions of $1$ in the Thue Morse sequence. Then
$$A+A=\mathbb{N}_0\setminus\{2,4,\,2^{2n+1}\!-1, n\ge 0\}.$$
$$B+B=\mathbb{N}_0\setminus\{2^{2n+1}\!-1, n\ge 0\}.$$
$$A+B=\mathbb{N}_0\setminus\{2^{2n}\!-1, n\ge 0\}.$$
\end{theorem}

 \medskip

 \emph{Proof:} We only prove the $A+A$ result. The other two can be proved in a similar way.

 According to Theorem \ref{th:diag},   $n\in A+A$ if and only if $d_{n-1}$ contains a symbol $(0,0)=a$. Note that we had to shift the diagonal words, as they are redefined for this section.


  Let $\sigma$ be the morphism from Proposition \ref{prop:diag} in $abcd$-coding,  $\sigma(a)=bc, \,\sigma(b)=ad, \,\sigma(c)=da, \,\sigma(d)=cb$. Then $\sigma^2$ is given by
 $$\sigma^2(a)=adda, \,\sigma^2(b)=bccb, \,\sigma^2(c)=cbbc, \,\sigma^2(d)=daad.$$
 It follows from this and Equation (\ref{eq:diagd15}) that $d(2^{2n+1})\in \{b,c\}^*,  \,d(2^{2n})\in \{a,d\}^*$, for $n\ge 0$.

 \medskip

This takes care of the $2^{2n+1}\!-1$ part of the theorem.  Obviously 2 and 4 are not in $A+A$. It remains to prove that all other numbers are in $A+A$, which is equivalent to proving that a symbol $a$  occurs in all  diagonal words  $d_n$ with $n>5$ and $n$ not equal to a power of 2.
It looks attractive to use Proposition \ref{prop:diag} to accomplish this. Indeed, let $T$ be the set consisting of the twelve 2-blocks $ab, ac,\dots,dc$ with two different symbols. Then one may check that all elements of $T$ occur in the length 4 blocks from $\sigma(T)$, and also in the length 3 blocks from $\beta(\sigma(T))$. This implies that a diagonal word $d_n$ in which all elements from $T$ occur, propagates this property to $d_{2n}$ and $d_{2n-1}$. However, there are many diagonal words in which not all blocks from $T$ occur, in fact all $n$ of the form $n=2^N+2^K$ for some non-negative integers $K$ and $N$. This makes an induction proof very complex, and so we will give a completely different proof.

\medskip

\begin{figure}[h!]%
\begin{center}
 \includegraphics[height=4cm]{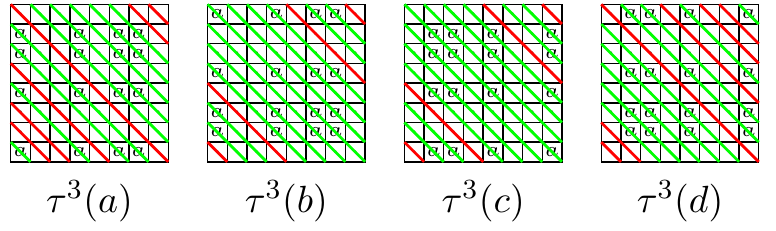}
 \vspace*{-1cm}
 \end{center}
 \caption{\footnotesize Diagonals in $\tau^3$-squares.}\label{Fig:ttp3}
\end{figure}

In Figure \ref{Fig:ttp3} we depict the  $\tau^3$-squares. Here the red lines indicate diagonal words \emph{without}   the symbol $a$, and the  green lines indicate diagonal words  \emph{with} a symbol $a$.

\smallskip

We also indicated (parts of) the diagonal words above the main diagonal. These are denoted by  $d_1^+,\dots,d_7^+$, starting from the main diagonal. For instance the block $\tau^3(a)$ has the two red diagonal lines $d_6^+$ and $d_7^+$.

\begin{figure}[h!]%
\begin{center}
 \includegraphics[height=4.8cm]{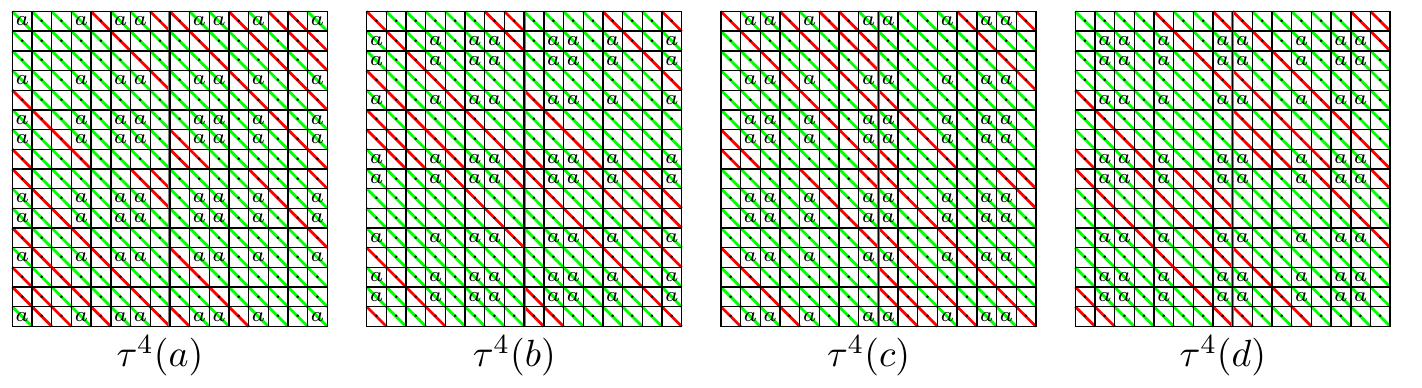}
 \vspace*{-1cm}
 \end{center}
 \caption{\footnotesize Diagonals in $\tau^4$-squares inherited from $\tau^3$-blocks.}\label{Fig:ttp4}
\end{figure}

\smallskip

Figure \ref{Fig:ttp4} gives the red and green line structure for the $\tau^4$-blocks inherited from the $\tau^3$-blocks.

\begin{figure}[h!]%
\begin{center}
\includegraphics[height=4.8cm]{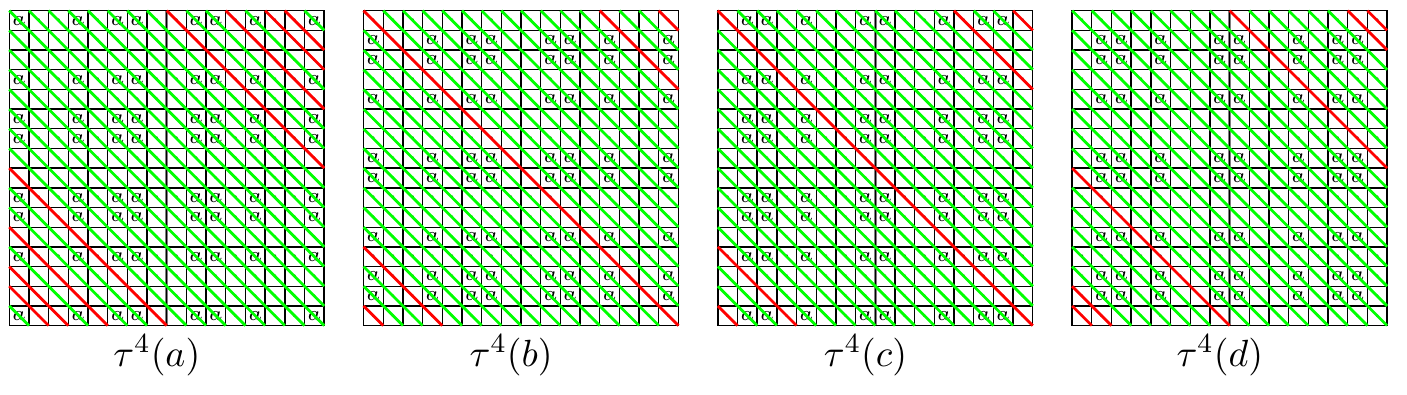}
 \vspace*{-1cm}
 \end{center}
 \caption{\footnotesize Final color assignment to the diagonals in $\tau^4$-blocks.}\label{Fig:ttp4-final}
\end{figure}

Observe that any diagonal that is partly red, partly green, actually should be a
green diagonal. This leads to Figure \ref{Fig:ttp4-final}.

\medskip

We finish the proof with induction. Starting from $n=4$, the blocks $\tau^n(a), \tau^n(b), \tau^n(c), \tau^n(d)$ have the following properties:

1) There are some red diagonals $d_i$ with $i\in \{1,2,3,4,5\}$,

2) There are some red diagonals $d^+_i$ with $i\in \{2^n-5,\dots,2^n-1\}$,

3) There are some red diagonals $d_i$ and $d^+_i$ with $i\in \{2^N, N\ge 2 \}$,

4) All other diagonals $d_i$ and $d^+_i$ are green.

\smallskip

One checks that these four properties hold for $n=4$, and then using the induction hypothesis makes the step from $n$ to $n+1$ in the same way as the step from $n=3$ to $n=4$ is made as in Figures \ref{Fig:ttp3}, \ref{Fig:ttp4} and \ref{Fig:ttp4-final}. \hfill $\Box$.


\section{The von Neumann word} \label{sec:vonNeu}

The von Neumann word $u$ was introduced in the paper \cite{All-Ber-Shallit}. One has $u=u_0u_1\dots=1101100110110\dots$, fixed point of the substitution \,$\nu:\quad 0\mapsto 0,\; 1\mapsto 110.$

The 2D von Neumann substitution $\psi=\nu\times\nu$ is  given by

\begin{tikzpicture}[scale=.4,rounded corners=0]
\def\bricka(#1){\path (#1)+(-1/2,-1/2) coordinate (shP0);
                \draw [fill=green!30!white, draw=black]  (shP0) rectangle +(1,1); \node[circle, inner sep=5.4pt] at (#1)  {\tiny $a$}}
\def\brickb(#1){\path (#1)+(-1/2,-1/2) coordinate (shP0);
                \draw [fill=yellow!30!white, draw=black]  (shP0) rectangle +(1,1); \node[circle, inner sep=5.4pt] at (#1)  {\tiny $b$}}
\def\brickc(#1){\path (#1)+(-1/2,-1/2) coordinate (shP0);
                \draw [fill=blue!30!white, draw=black]  (shP0) rectangle +(1,1); \node[circle, inner sep=5.4pt] at (#1)  {\tiny $c$}}
\def\brickd(#1){\path (#1)+(-1/2,-1/2) coordinate (shP0);
                \draw [fill=red!30!white, draw=black]  (shP0) rectangle +(1,1); \node[circle, inner sep=5.4pt] at (#1)  {\tiny $d$}}
\def\theta(#1,#2){\path (#2) coordinate (P0);
\path (P0)++(0,0) coordinate (P00); \path (P0)++(0,1) coordinate (P01); \path (P0)++(0,2) coordinate (P02);
\path (P0)++(1,0) coordinate (P10); \path (P0)++(1,1) coordinate (P11); \path (P0)++(1,2) coordinate (P12);
\path (P0)++(2,0) coordinate (P20); \path (P0)++(2,1) coordinate (P21); \path (P0)++(2,2) coordinate (P22);
\ifthenelse{#1=00}{\bricka(P00) }{};                                 
\ifthenelse{#1=01}{\brickb(P00); \brickb(P01); \bricka(P02) }{};     
\ifthenelse{#1=10}{\brickc(P00); \brickc(P10); \bricka(P20) }{};     
\ifthenelse{#1=11}{\brickd(P00); \brickd(P01); \brickc(P02); \brickd(P10); \brickd(P11); \brickc(P12); \brickb(P20); \brickb(P21); \bricka(P22) }{};
}; 
\path (0,-6) coordinate (a);   \bricka(a);  \node at (2,-6) {$\mapsto$};  \path (4,-6) coordinate (t00);   \theta(00,t00);
\path (7,-6) coordinate (b);   \brickb(b);  \node at (9,-6) {$\mapsto$};  \path (11,-6) coordinate (t01);  \theta(01,t01);
\path (14,-6) coordinate (c);  \brickc(c);  \node at (16,-6) {$\mapsto$}; \path (18,-6) coordinate (t10);  \theta(10,t10);
\path (23,-6) coordinate (d);  \brickd(d);  \node at (25,-6) {$\mapsto$}; \path (27,-6) coordinate (t11);  \theta(11,t11);
\end{tikzpicture}

Before we continue, we formulate a simple lemma on the language $L_\nu$ of the substitution $\nu$, i.e., on the collection of words that can occur in $u$.

\begin{lemma}\label{lem:words}
The words $010, 111$ and $101101$ do not occur in $L_\nu$.
\end{lemma}

\smallskip

Using Lemma \ref{lem:words} one simply derives the next lemma.

\begin{lemma}\label{lem:wordscases}
The word $11$ occurs uniquely as suffix of $011$. The word $101$ occurs uniquely as suffix of $001101$. The word $10110$ occurs uniquely as suffix of $0110110$.
\end{lemma}

\begin{figure}[h!]%
\begin{center}
\includegraphics[height=6.5cm]{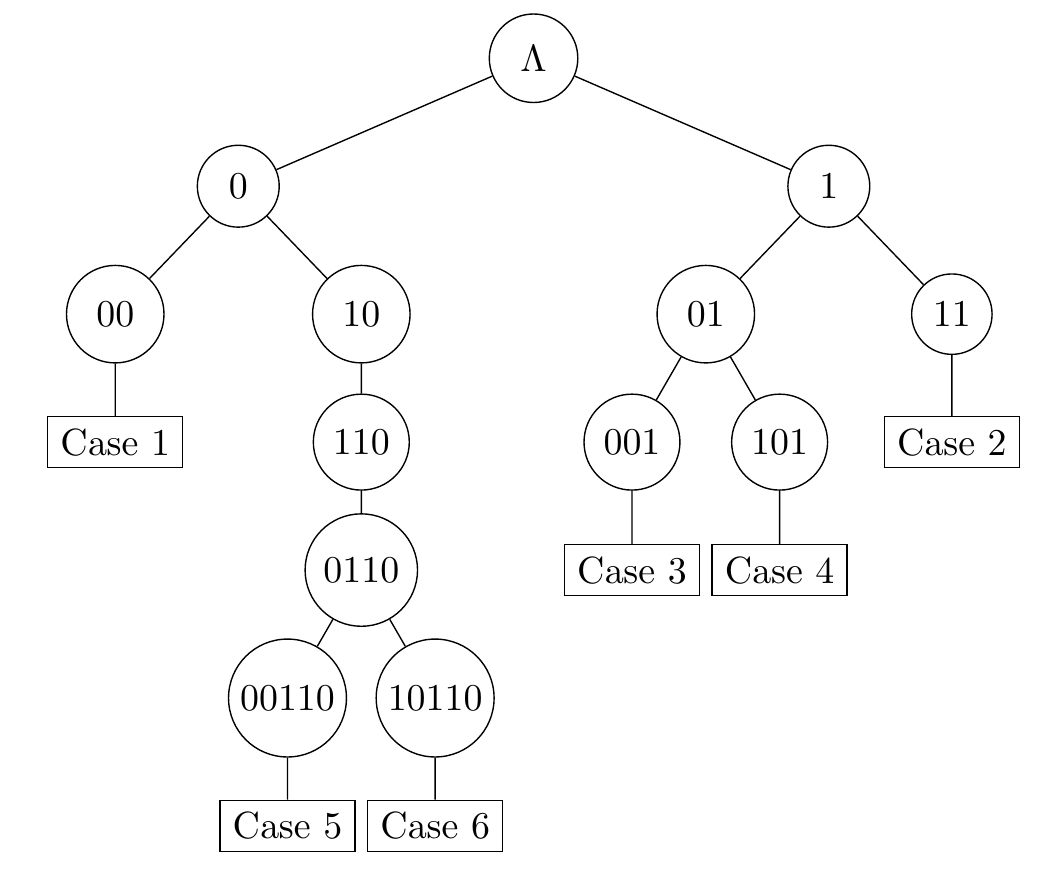}
 \vspace*{-.3cm}
 \end{center}
 \caption{\footnotesize The prefix tree for the von Neumann words.}\label{Fig:prefixtree}
\end{figure}
 \begin{figure}[h!]%
\begin{center}
\includegraphics[height=4.2cm]{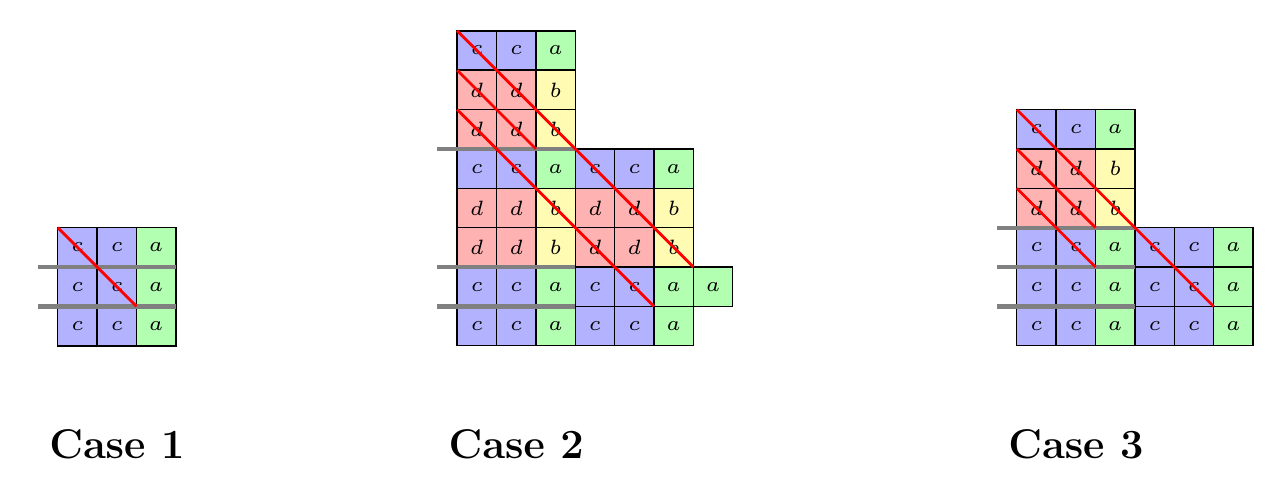}
 \vspace*{-.8cm}
 \end{center}
 \caption{\footnotesize Diagonals for 2D von Neumann words in Case 1,2,3.}\label{Fig:case123}
\end{figure}

We split the possible occurrences of words of $L_\nu$ in the von Neumann word $u$ according to a tree of prefixes, dividing these in six possible cases. This prefix tree is given in Figure \ref{Fig:prefixtree}, where $\Lambda$ is the empty word. Note that we used Lemma \ref{lem:words} again to determine the offspring of the node labeled 10.

\begin{theorem}\label{th:vN}
Let $A=3,6,7,10,13,\dots$ give the positions of $0$'s in the von Neumann sequence, and $B=1,2,4,5,8,9\dots$ the positions of $1$ in the von Neumann sequence. Then
$$A+A=\mathbb{N}\setminus\{2,3,4,5,7,8,11,15\}.$$
$$B+B=\mathbb{N}\setminus\{1\}.$$
\end{theorem}

\medskip

\noindent \emph{Proof:} The proof of $B+B=\mathbb{N}\setminus\{1\}$ is almost trivial.

 The proof for $A+A$ is based on Theorem \ref{th:diag}, with the trivial change that we consider the fixed point starting with 1. The proof runs as the proof for the Fibonacci word, except that here we need the 15 left most columns.

 First one checks easily that the numbers $2,3,4,5,7,8,11,15$ are not in $A+A$.

  To handle the diagonal words $d_n$ for $n\ge 16$, we split into the six cases we introduced in Figure \ref{Fig:prefixtree}. In Case 2, 3 and 4 we consider the three diagonals starting at the left border in the blocks $\psi(d)$, in the remaining four the single diagonal starting in $\psi(c)$. See Figure \ref{Fig:case123} for the three cases 00, 11 and 001. This means that in Case 1 the  top two blocks are $\psi(c)\psi(c)$, in Case 2 they are $\psi(d)\psi(d)$, and in Case 3 they are $\psi(c)\psi(c)\psi(d)$.

  \begin{figure}[h!]%
\begin{center}
\includegraphics[height=6.0cm]{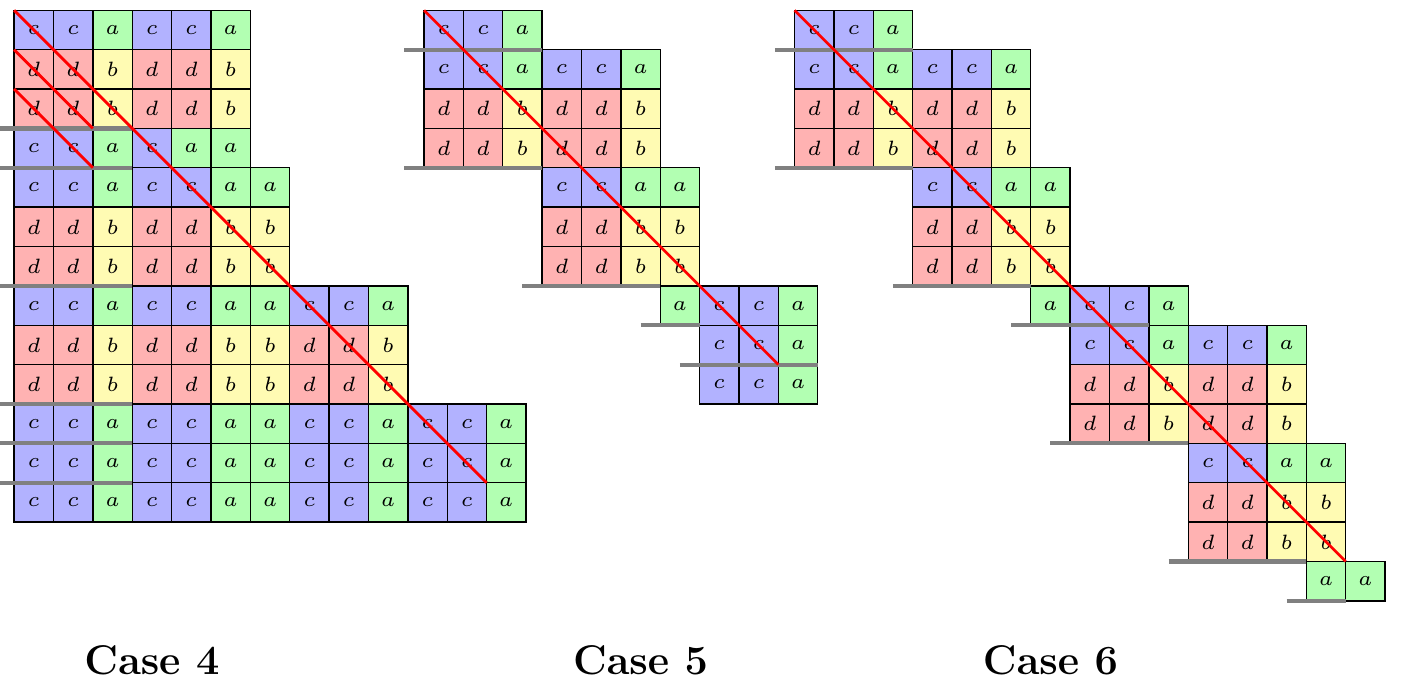}
 \vspace*{-.8cm}
 \end{center}
 \caption{\footnotesize Diagonals for 2D von Neumann words in  Case 4,5,6.}\label{Fig:case456}
\end{figure}

In Case 2  a block $\psi(c)$, must precede the top blocks $\psi(d)\psi(d)$, by Lemma 2. In a similar way there are forced blocks in Case 4 and Case 6.

The other blocks are then inserted following the von Neumannn words $ccaccaa\dots$ or $ddbddbb\dots$. In all cases, except the last, a row $ccaccaacca\dots$ has been added at the bottom. This is allowed because both $\psi(c)$ and $\psi(d)$ have the letter $c$ as prefix.

The red diagonals are those parts of the three diagonals starting in the top block $\psi(d)$, respectively the single diagonal starting in $\psi(c)$,  that do not cross a square with label $a$, i.e., the red diagonal ends just before the first $a$-square.
See Figure \ref{Fig:case456} for the three cases remaining cases.

The fact that in all six cases a square with label $a$ is encountered finishes the proof. \hfill $\Box$

\section{The sum of squares} \label{sec:sq}

Here we give a classical example of a sumset. Let $A=\{n^2: n\ge 1\}$. It is well known (see, e.g., \cite{Cobham}), that the characteristic function of $A$, as a word, is a letter to letter substitution $\lambda$ of the fixed point with prefix 0 of the morphism
$$ 0\mapsto 01, \quad  1\mapsto 221, \quad  2\mapsto 2.$$
The letter to letter map is given by $\lambda(0)=0, \;\lambda(1)=1,\;\lambda(2)=0$.

\medskip

The corresponding 2D substitution is the morphism $\mu$ given by

\bigskip

\begin{tikzpicture}[scale=.4,rounded corners=0]
\def\brickZZ(#1){\path (#1)+(-1/2,-1/2) coordinate (shP0);
                \draw [fill=green!40!white, draw=gray,thick]  (shP0) rectangle +(1,1); \node[circle, inner sep=5.4pt] at (#1)  {\tiny $00$}};
\def\brickZO(#1){\path (#1)+(-1/2,-1/2) coordinate (shP0);
                \draw [fill=yellow!60!white, draw=gray,thick]  (shP0) rectangle +(1,1); \node[circle, inner sep=5.4pt] at (#1)  {\tiny $01$}} ;
\def\brickZT(#1){\path (#1)+(-1/2,-1/2) coordinate (shP0);
                \draw [fill=yellow!30!white, draw=gray,thick]  (shP0) rectangle +(1,1); \node[circle, inner sep=5.4pt] at (#1)  {\tiny $02$}} ;
\def\brickOZ(#1){\path (#1)+(-1/2,-1/2) coordinate (shP0);
                \draw [fill=yellow!60!white, draw=gray,thick]  (shP0) rectangle +(1,1); \node[circle, inner sep=5.4pt] at (#1)  {\tiny $10$}} ;
\def\brickOO(#1){\path (#1)+(-1/2,-1/2) coordinate (shP0);
                \draw [fill=red!60!white, draw=gray,thick]  (shP0) rectangle +(1,1); \node[circle, inner sep=5.4pt] at (#1)  {\tiny $11$}};
\def\brickOT(#1){\path (#1)+(-1/2,-1/2) coordinate (shP0);
                \draw [fill=blue!30!white, draw=gray,thick]  (shP0) rectangle +(1,1); \node[circle, inner sep=5.4pt] at (#1)  {\tiny $12$}} ;
\def\brickTZ(#1){\path (#1)+(-1/2,-1/2) coordinate (shP0);
                \draw [fill=yellow!30!white, draw=gray,thick]  (shP0) rectangle +(1,1); \node[circle, inner sep=5.4pt] at (#1)  {\tiny $20$}} ;
\def\brickTO(#1){\path (#1)+(-1/2,-1/2) coordinate (shP0);
                \draw [fill=blue!30!white, draw=gray,thick]  (shP0) rectangle +(1,1); \node[circle, inner sep=5.4pt] at (#1)  {\tiny $21$}} ;
\def\brickTT(#1){\path (#1)+(-1/2,-1/2) coordinate (shP0);
                \draw [fill=blue!10!white, draw=gray,thick]  (shP0) rectangle +(1,1); \node[circle, inner sep=5.4pt] at (#1)  {\tiny $22$}} ;
\def\theta(#1,#2){\path (#2) coordinate (P0);
\path (P0)++(0,0) coordinate  (P00);
\path (P0)++(1,0) coordinate (P10); \path (P0)++(0,1) coordinate (P01); \path (P0)++(1,1) coordinate (P11); \path (P0)++(2,0) coordinate (P20); \path (P0)++(0,2) coordinate (P02); \path (P0)++(1,2) coordinate (P12); \path (P0)++(2,1) coordinate (P21); \path (P0)++(2,2) coordinate (P22);
\ifthenelse{#1=00}{\brickZZ(P00);  \brickZO(P01); \brickOZ(P10); \brickOO(P11)}{};
\ifthenelse{#1=01}{\brickZT(P00);  \brickZT(P01); \brickOT(P10); \brickOT(P11); \brickZO(P02); \brickOO(P12)}{};
\ifthenelse{#1=02}{\brickZT(P00);  \brickOT(P10) }{};
\ifthenelse{#1=10}{\brickTZ(P00);  \brickTO(P01); \brickTO(P11); \brickTZ(P10);  \brickOO(P21); \brickOZ(P20)}{};
\ifthenelse{#1=11}{\brickTT(P00);   \brickTT(P01);  \brickTO(P02);
                   \brickTT(P10);   \brickTT(P11);  \brickTO(P12);
                   \brickOT(P20);   \brickOT(P21);  \brickOO(P22) }{};
\ifthenelse{#1=12}{\brickTT(P00);  \brickTT(P10); \brickOT(P20)}{};
\ifthenelse{#1=20}{\brickTZ(P00);  \brickTO(P01); }{};
\ifthenelse{#1=21}{\brickTT(P00);  \brickTT(P01); \brickTO(P02) }{};
\ifthenelse{#1=22}{\brickTT(P00); }{};
}; 

\path (0,-4) coordinate (ZZ);   \brickZZ(ZZ);  \node at (2,-4) {$\mapsto$};  \path (4,-4) coordinate (t00);  \theta(00,t00);
\path (10,-4) coordinate (ZO);   \brickZO(ZO);  \node at (12,-4) {$\mapsto$};  \path (14,-4) coordinate (t01);  \theta(01,t01);
\path (20,-4) coordinate (ZT);  \brickZT(ZT);  \node at (22,-4) {$\mapsto$};  \path (24,-4) coordinate (t02);  \theta(02,t02);
\path (0,-8) coordinate (OZ);  \brickOZ(OZ);  \node at (2,-8) {$\mapsto$};  \path (4,-8) coordinate (t10); \theta(10,t10);
\path (10,-8) coordinate (OO);  \brickOO(OO);  \node at (12,-8) {$\mapsto$}; \path (14,-8) coordinate (t11); \theta(11,t11);
\path (20,-8) coordinate (OT);  \brickOT(OT);  \node at (22,-8) {$\mapsto$}; \path (24,-8) coordinate (t12); \theta(12,t12);
\path (0,-12) coordinate (TZ);  \brickTZ(TZ);  \node at (2,-12) {$\mapsto$}; \path (4,-12) coordinate (t20); \theta(20,t20);
\path (10,-12) coordinate (TO);  \brickTO(TO);  \node at (12,-12) {$\mapsto$}; \path (14,-12) coordinate (t21); \theta(21,t21);
\path (20,-12) coordinate (TT);  \brickTT(TT);  \node at (22,-12) {$\mapsto$}; \path (24,-12) coordinate (t22); \theta(22,t22);
\end{tikzpicture}

\medskip

The numbers $s=n^2+m^2$ that are a sum of two squares occur as sums of the indices $(n^2,m^2)$  of the squares with the symbols $11$ in the fixed point $S$ of $\mu$, limit of $\mu^n(00)$ as $n\rightarrow \infty$.

\smallskip

At first sight it is surprising that this set has a fractal structure, but due
to the fact that the morphism $\mu$ is not primitive (i.e., the incidence matrix of the substitution is reducible), there is no exponential scaling structure in the 2D word $S$, but rather a polynomial one. The latter simply amounts to the recursion $(n+1)^2=n^2+(2n+1)$.
This recursion is clearly visible in  Figure \ref{Fig:sq}, which displays the two-dimensional word  $\mu^4(00)$.

\medskip
\begin{figure}[h!]%
\begin{center}
\begin{tikzpicture}[scale=.4,rounded corners=0]
\def\brickZZ(#1){\path (#1)+(-1/2,-1/2) coordinate (shP0);
                \draw [fill=green!40!white, draw=gray,thick]  (shP0) rectangle +(1,1); \node[circle, inner sep=5.4pt] at (#1)  {\tiny $00$}};
\def\brickZO(#1){\path (#1)+(-1/2,-1/2) coordinate (shP0);
                \draw [fill=yellow!60!white, draw=gray,thick]  (shP0) rectangle +(1,1); \node[circle, inner sep=5.4pt] at (#1)  {\tiny $01$}} ;
\def\brickZT(#1){\path (#1)+(-1/2,-1/2) coordinate (shP0);
                \draw [fill=yellow!30!white, draw=gray,thick]  (shP0) rectangle +(1,1); \node[circle, inner sep=5.4pt] at (#1)  {\tiny $02$}} ;
\def\brickOZ(#1){\path (#1)+(-1/2,-1/2) coordinate (shP0);
                \draw [fill=yellow!60!white, draw=gray,thick]  (shP0) rectangle +(1,1); \node[circle, inner sep=5.4pt] at (#1)  {\tiny $10$}} ;
\def\brickOO(#1){\path (#1)+(-1/2,-1/2) coordinate (shP0);
                \draw [fill=red!60!white, draw=gray,thick]  (shP0) rectangle +(1,1); \node[circle, inner sep=5.4pt] at (#1)  {\tiny $11$}};
\def\brickOT(#1){\path (#1)+(-1/2,-1/2) coordinate (shP0);
                \draw [fill=blue!30!white, draw=gray,thick]  (shP0) rectangle +(1,1); \node[circle, inner sep=5.4pt] at (#1)  {\tiny $12$}} ;
\def\brickTZ(#1){\path (#1)+(-1/2,-1/2) coordinate (shP0);
                \draw [fill=yellow!30!white, draw=gray,thick]  (shP0) rectangle +(1,1); \node[circle, inner sep=5.4pt] at (#1)  {\tiny $20$}} ;
\def\brickTO(#1){\path (#1)+(-1/2,-1/2) coordinate (shP0);
                \draw [fill=blue!30!white, draw=gray,thick]  (shP0) rectangle +(1,1); \node[circle, inner sep=5.4pt] at (#1)  {\tiny $21$}} ;
\def\brickTT(#1){\path (#1)+(-1/2,-1/2) coordinate (shP0);
                \draw [fill=blue!10!white, draw=gray,thick]  (shP0) rectangle +(1,1); \node[circle, inner sep=5.4pt] at (#1)  {\tiny $22$}} ;
\def\theta(#1,#2){\path (#2) coordinate (P0);
\path (P0)++(0,0) coordinate  (P00);
\path (P0)++(1,0) coordinate (P10); \path (P0)++(0,1) coordinate (P01); \path (P0)++(1,1) coordinate (P11); \path (P0)++(2,0) coordinate (P20); \path (P0)++(0,2) coordinate (P02); \path (P0)++(1,2) coordinate (P12); \path (P0)++(2,1) coordinate (P21); \path (P0)++(2,2) coordinate (P22);
\ifthenelse{#1=00}{\brickZZ(P00);  \brickZO(P01); \brickOZ(P10); \brickOO(P11)}{};
\ifthenelse{#1=01}{\brickZT(P00);  \brickZT(P01); \brickOT(P10); \brickOT(P11); \brickZO(P02); \brickOO(P12)}{};
\ifthenelse{#1=02}{\brickZT(P00);  \brickOT(P10) }{};
\ifthenelse{#1=10}{\brickTZ(P00);  \brickTO(P01); \brickTO(P11); \brickTZ(P10);  \brickOO(P21); \brickOZ(P20)}{};
\ifthenelse{#1=11}{\brickTT(P00);   \brickTT(P01);  \brickTO(P02);
                   \brickTT(P10);   \brickTT(P11);  \brickTO(P12);
                   \brickOT(P20);   \brickOT(P21);  \brickOO(P22) }{};
\ifthenelse{#1=12}{\brickTT(P00);  \brickTT(P10); \brickOT(P20)}{};
\ifthenelse{#1=20}{\brickTZ(P00);  \brickTO(P01); }{};
\ifthenelse{#1=21}{\brickTT(P00);  \brickTT(P01); \brickTO(P02) }{};
\ifthenelse{#1=22}{\brickTT(P00); }{};
}; 
\def\thetatwo(#1,#2){\path (#2) coordinate (Q0);
\path (Q0)++(0,0) coordinate (Q00); \path (Q0)++(2,0) coordinate (Q10); \path (Q0)++(4,0) coordinate (Q20);
\path (Q0)++(0,1) coordinate (Q01); \path (Q0)++(2,1) coordinate (Q11); \path (Q0)++(4,4) coordinate (Q21);
\path (Q0)++(0,2) coordinate (Q02); \path (Q0)++(2,2) coordinate (Q12); \path (Q0)++(2,2) coordinate (Q22);
\ifthenelse{#1=00}{ \path (Q0)++(0,2) coordinate (Q01); \path (Q0)++(2,2) coordinate (Q11);
    \theta(00,Q00);  \theta(01,Q01); \theta(10,Q10); \theta(11,Q11)}{};
\ifthenelse{#1=01}{\theta(02,Q00);  \theta(02,Q01); \theta(12,Q10); \theta(12,Q11); \theta(01,Q02); \theta(11,Q12)}{};
\ifthenelse{#1=02}{\theta(02,Q00);  \theta(12,Q10) }{};
\ifthenelse{#1=10}{\path (Q0)++(1,0) coordinate (Q10); \path (Q0)++(2,0) coordinate (Q20);
   \path (Q0)++(0,2) coordinate (Q01); \path (Q0)++(1,2) coordinate (Q11); \path (Q0)++(2,2) coordinate (Q21);
   \theta(20,Q00);  \theta(20,Q10); \theta(21,Q11); \theta(21,Q01);  \theta(11,Q21); \theta(10,Q20)}{};
\ifthenelse{#1=11}{\path (Q0)++(1,0) coordinate (Q10); \path (Q0)++(2,0) coordinate (Q20);
   \path (Q0)++(0,1) coordinate (Q01); \path (Q0)++(1,1) coordinate (Q11); \path (Q0)++(2,1) coordinate (Q21);
   \path (Q0)++(0,2) coordinate (Q02); \path (Q0)++(1,2) coordinate (Q12); \path (Q0)++(2,2) coordinate (Q22);
                   \theta(22,Q00);  \theta(22,Q01);  \theta(21,Q02);
                   \theta(22,Q10);  \theta(22,Q11);  \theta(21,Q12);
                   \theta(12,Q20);  \theta(12,Q21);  \theta(11,Q22) }{};
\ifthenelse{#1=12}{\path (Q0)++(1,0) coordinate (Q10); \path (Q0)++(2,0) coordinate (Q20);
      \theta(22,Q00);  \theta(22,Q10); \theta(12,Q20)}{};
\ifthenelse{#1=20}{ \path (Q0)++(0,2) coordinate (Q01); \theta(20,Q00);  \theta(21,Q01); }{};
\ifthenelse{#1=21}{\theta(22,Q00);  \theta(22,Q01); \theta(21,Q02) }{};
\ifthenelse{#1=22}{\theta(22,Q00); }{}
}; 
\def\thetathree(#1,#2){\path (#2) coordinate (R0);
\path (R0)++(0,0) coordinate  (R00); \path (R0)++(5,0) coordinate (R10); \path (R0)++(0,5) coordinate (R01); \path (R0)++(5,5) coordinate (R11);
\path (R0)++(0,5) coordinate (R02); \path (R0)++(5,5) coordinate (R12); \path (R0)++(5,5) coordinate (R22);
\ifthenelse{#1=00}{\thetatwo(00,R00);  \thetatwo(01,R01); \thetatwo(10,R10); \thetatwo(11,R11)}{};
\ifthenelse{#1=01}{\path (R0)++(5,0) coordinate (R10); \path (R0)++(0,1) coordinate (R01);
   \path (R0)++(5,1) coordinate (R11); \path (R0)++(0,2) coordinate (R02);  \path (R0)++(5,2) coordinate (R12);
    \thetatwo(02,R00);  \thetatwo(02,R01); \thetatwo(12,R10); \thetatwo(12,R11); \thetatwo(01,R02); \thetatwo(11,R12)}{};
\ifthenelse{#1=02}{\thetatwo(02,R00);  \thetatwo(12,R10) }{};
\ifthenelse{#1=10}{\path (R0)++(1,0) coordinate (R10); \path (R0)++(2,0) coordinate (R20);
   \path (R0)++(0,5) coordinate (R01); \path (R0)++(1,5) coordinate (R11); \path (R0)++(2,5) coordinate (R21);
   \thetatwo(20,R00);  \thetatwo(20,R10); \thetatwo(21,R11); \thetatwo(21,R01);  \thetatwo(11,R21); \thetatwo(10,R20)}{};
\ifthenelse{#1=11}{\path (R0)++(1,0) coordinate (R10); \path (R0)++(2,0) coordinate (R20);
   \path (R0)++(0,1) coordinate (R01); \path (R0)++(1,1) coordinate (R11); \path (R0)++(2,1) coordinate (R21);
   \path (R0)++(0,2) coordinate (R02); \path (R0)++(1,2) coordinate (R12); \path (R0)++(2,2) coordinate (R22);
                   \thetatwo(22,R00);  \thetatwo(22,R01);  \thetatwo(21,R02);
                   \thetatwo(22,R10);  \thetatwo(22,R11);  \thetatwo(21,R12);
                   \thetatwo(12,R20);  \thetatwo(12,R21);  \thetatwo(11,R22) }{};
\ifthenelse{#1=12}{\path (R0)++(1,0) coordinate (R10); \path (R0)++(2,0) coordinate (R20);
      \thetatwo(22,R00);  \thetatwo(22,R10); \thetatwo(12,R20)}{};
\ifthenelse{#1=20}{\thetatwo(20,R00); \thetatwo(21,R01); }{};
\ifthenelse{#1=21}{\path (R0)++(0,1) coordinate (R01); \path (R0)++(0,2) coordinate (R02);
    \thetatwo(22,R00);  \thetatwo(22,R01); \thetatwo(21,R02) }{};
\ifthenelse{#1=22}{\thetatwo(22,R00); }{}
}; 
\def\thetafour(#1,#2){\path (#2) coordinate (S0);
\path (S0)++(0,0) coordinate  (S00); \path (S0)++(10,0) coordinate (S10); \path (S0)++(0,10) coordinate (S01); \path (S0)++(10,10) coordinate (S11);
\ifthenelse{#1=00}{\thetathree(00,S00);  \thetathree(01,S01); \thetathree(10,S10); \thetathree(11,S11) }{};
}; 
\thetafour(00,(0,31));
\vspace*{-1cm}
\end{tikzpicture}
 \end{center}
 \caption{\footnotesize The 2D word $\mu^4(00)$.}\label{Fig:sq}
\end{figure}

\section{Acknowledgement}
I am grateful to Jean-Paul Allouche and Jeffrey Shallit for their useful comments on this paper.

\noindent {\footnotesize AMS Classification Numbers: 11B13, 68R10.}

\end{document}